\documentclass[12pt]{amsart}
\usepackage{amsmath}
\usepackage{amssymb}

\newcommand{\cz}{{\mathbb C}}
\newcommand{\fz}{{\mathbb F}}
\newcommand{\pz}{{\mathbb P}}
\newcommand{\qz}{{\mathbb Q}}
\newcommand{\zz}{{\mathbb Z}}
\newcommand{\lra}{{\longrightarrow}}

\newcommand{\Frob}{\operatorname{Frob}} 
\newcommand{\tr}{\operatorname{tr}}     

\newtheorem{Theorem}{Theorem}[section]
\newtheorem{Question}[Theorem]{Question}

\begin{document}

\title{The mirror quintic as a quintic}

\author{Christian Meyer}
\thanks{Partially supported by DFG Schwerpunktprogramm 1094
(Globale Methoden in der komplexen Geometrie). The author also thanks
D. van Straten for helpful discussions.}

\address{Fachbereich Mathematik und Informatik\\ Johannes
  Gutenberg-Uni\-ver\-sit\"at\\
  Staudingerweg 9\\D--55099 Mainz\\Germany}
\email{cm@mathematik.uni-mainz.de}

\maketitle

\section{Introduction}

The general quintic hypersurface in $\pz^4$ is the most famous example of a Calabi--Yau
threefold for which mirror symmetry has been investigated in detail. There is a
description of the mirror as a hypersurface in a certain weighted projective space.
In this note we present a model for the mirror which is
again (the resolution of) a quintic hypersurface in $\pz^4$. We also deal with the
special members in the respective families. They lead to rigid Calabi--Yau threefolds
with interesting arithmetical properties.

In the last section we try to find a similarly nice model for the mirror of the
complete intersection of two cubics in $\pz^5$. We also formulate a general
question about mirror models.

Some of the results of this note are contained in the author's thesis (\cite{Meyer})
which can also be used as an introduction to modularity of Calabi--Yau threefolds.

\section{The quintic}

Let $X_{\mu} \subset \pz^4$ be the quintic threefold defined by the equation
\[
x_0^5 + x_1^5 + x_2^5 + x_3^5 + x_4^5 - 5 \mu x_0 x_1 x_2 x_3 x_4 = 0.
\]

If $\mu$ is general (i.e., no 5-th root of unity and not 0 or $\infty$) then $X_{\mu}$
is smooth and so a Calabi--Yau threefold. Its invariants are
\[
\chi(X_{\mu}) = -200, \qquad h^{1,1}(X_{\mu}) = 1, \qquad h^{2,1}(X_{\mu}) = 101.
\]

On $X_{\mu}$ there is an action of the group $G\simeq (\zz/5\zz)^3$ generated by
the coordinate transformations
\[
(x_0:x_1:x_2:x_3:x_4:x_5) \mapsto
(x_0:x_1 \cdot \xi_5^{\lambda_1}:x_2 \cdot \xi_5^{\lambda_2}:
x_3 \cdot \xi_5^{\lambda_3}:x_4 \cdot \xi_5^{\lambda_4})
\]
with $\lambda_1, \lambda_2, \lambda_3, \lambda_4 \in \zz/5\zz$,
$\lambda_1 + \lambda_2 + \lambda_3 + \lambda_4 \equiv 0 \mod 5$
and $\xi_5$ a fixed primitive 5-th root of unity.

The mirror partner of $X_{\mu}$ can be described as a resolution of the quotient $X_{\mu}/G$.
In fact this example was the first for which mirror symmetry was used to make enumerative
predictions (Candelas, de la Ossa, Green and Parkes in \cite{CanOssGP}), and it is still
by far the most prominent example in the mirror symmetry business.
There is a lot of literature; a good starting point
is the book of Cox and Katz (\cite{CoxKatz}, chapter 2 and section 4.2).
Morrison (\cite{Morrison}) describes in detail the
resolution of singularities of the mirror, and Candelas, de la Ossa and Rodriguez-Villegas
(\cite{CanOssVil1}, \cite{CanOssVil2}) discuss the zeta functions.

If $\mu$ is a 5-th root of unity then $X_{\mu}$ has 125 ordinary double points
($A_1$ singularities, also called {\em nodes}) as only singularities,
namely the points on the orbit of the point $(1:\mu:\mu:\mu:\mu)$ under the action
of the group $G$. The threefold $X:=X_1$ is defined over $\qz$ and especially interesting
for arithmetical purposes. Schoen (\cite{Schoen}) discussed the {\em modularity} of $X$,
i.e., he determined the connection between the arithmetic of $X$ over finite fields
and a certain weight four modular form. We will briefly recall his results.

Let $\tilde{X}$ be a small resolution of $X$. Then $\tilde{X}$ has Euler
characteristic
\[
\chi(\tilde{X}) = -200 + 2\cdot 125 = 50
\]
and Hodge numbers
\[
h^{1,1}(\tilde{X}) = 25, \qquad h^{2,1}(\tilde{X}) = 0.
\]

The defect of $X$ is $d(X)=h^{1,1}(\tilde{X}) - 1=24 > 0$. Since the group $G$ acts
transitively on the set of nodes of $X$ there exist projective small resolutions
(cf. \cite{Werner}, chapter IV).
Equivalently this may be deduced from the existence of smooth quadric surfaces on $X$
containing all the nodes (cf. \cite{Schoen}). For instance, the node $(1:1:1:1:1)$ is contained
in the smooth quadric surface $Q$ given by the equations
\begin{align*}
& x_0 + \xi_5 x_1 + \xi_5^2 x_2 + \xi_5^3 x_3 + \xi_5^4 x_4 =\\
& \qquad x_0 x_1 + \xi_5 x_0 x_2 + \xi_5^2 x_0 x_3 + \xi_5^3 x_0 x_4 + \xi_5^2 x_1 x_2\\
& \qquad + \xi_5^3 x_1 x_3 + \xi_5^4 x_1 x_4 + \xi_5^4 x_2 x_3 + x_2 x_4 + \xi_5 x_3 x_4 = 0 
\end{align*}
where $\xi_5$ is a primitive 5-th root of unity.

Schoen (\cite{Schoen}) proves that for all primes $p$ of good reduction for $X$
(i.e., $p\neq 5$) we have
\begin{align*}
& \tr(\Frob_p^*|H^{3}_{\text{\'et}}(\tilde{X},\qz_{\ell})) = a_p\\
& \qquad = \begin{cases}
p^3 + 25 p^2 - 100p + 1 - \#X_p, & p \equiv 1 \mod 5,\\
p^3 + \phantom{01}p^2 \hspace*{38pt} + 1 - \#X_p, & p \equiv 4 \mod 5,\\
p^3 + \phantom{01}p^2 + \phantom{00}2p + 1 - \#X_p, & p \equiv 2,3 \mod 5.\\
\end{cases}
\end{align*}

Here $\Frob_p^*$ denotes the map on $l$-adic cohomology which is induced by the geometric
Frobenius, $\#X_p$ denotes the number of points on $X$ over the finite field $\fz_p$,
and $a_p$ are the coefficients of a certain weight four newform for $\Gamma_0(25)$.

\section{The mirror}

Let $Y_{\mu} \subset \pz^4$ be the quintic threefold defined by the equation
\[
( x_0 + x_1 + x_2 + x_3 + x_4 )^5 - (5\mu)^5 x_0 x_1 x_2 x_3 x_4 = 0.
\]

For general $\mu$ (i.e., $\mu$ is no 5-th root of unity and not 0 or $\infty$)
the singular locus of $Y_{\mu}$ consists of the 10 lines given by
\[
x_i = x_j = x_k + x_l + x_m = 0
\]
where $\{i,j,k,l,m\} = \{0,1,2,3,4\}$. If $\mu$ is a 5-th of unity then there is an
additional singularity at the point $(1:1:1:1:1)$. This is an ordinary node.

There is a rational dominant map $X_{\mu}\lra Y_{\mu}$ induced by
\[
\phi : \pz^4 \lra \pz^4, \quad (x_0:x_1:x_2:x_3:x_4) \mapsto (x_0^5:x_1^5:x_2^5:x_3^5:x_4^5).
\]

The map is generically $125:1$. The degree reduces to $25:1$ on the singular lines
and to $5:1$ on the 10 intersection points of three lines (i.e., the points on the
orbit of $(0:0:0:1:-1)$ under permutation of coordinates). Let $A$ denote the union
of the 10 singular lines and let $B$ denote the set containing the 10 intersection points.
We can now relate the Euler characteristic of a general $Y_{\mu}$ to that of $X_{\mu}$:
\begin{align*}
-200 = \chi(X_{\mu}) & = 125\cdot \chi(Y_{\mu}\setminus A) + 25\cdot \chi(A\setminus B) + 5\cdot \chi(B)\\
& = 125\cdot \chi(Y_{\mu}\setminus A) + 25\cdot 10\cdot (2-3) + 5\cdot 10\\
& = 125\cdot \chi(Y_{\mu}\setminus A) - 200,
\end{align*}
thus $\chi(Y_{\mu}\setminus A) = 0$ and
\[
\chi(Y_{\mu}) = \chi(Y_{\mu}\setminus A) + \chi(A) = 0 + 10\cdot 2 - 2\cdot 10 = 0.
\]

If $\mu$ is a 5-th root of unity then we set $Y:=Y_{\mu}=Y_1$ and we find
\begin{align*}
-75 = \chi(X) & = 125\cdot \chi(Y\setminus A) + 25\cdot \chi(A\setminus B) + 5\cdot \chi(B)\\
& = 125\cdot \chi(Y\setminus A) + 25\cdot 10\cdot (2-3) + 5\cdot 10\\
& = 125\cdot \chi(Y\setminus A) - 200,
\end{align*}
thus $\chi(Y\setminus A) = 1$ and
\[
\chi(Y) = \chi(Y\setminus A) + \chi(A) = 1 + 10\cdot 2 - 2\cdot 10 = 1.
\]

The map $\phi$ exactly divides out the action of the group $G$ on $X_{\mu}$.
Thus the quintic $Y_{\mu}$ is a model for the mirror $X_{\mu}/G$ of $X_{\mu}$.
The resolution of singularities of $X_{\mu}/G$ has been discussed in \cite{Morrison}.
The singular lines are lines of $A_4$ singularities. The 10 intersection points of
singular lines look locally like the quotient $\cz^3/H$ where the group
$H\cong\{(\xi_1,\xi_2,\xi_3),\, \xi_1^5 = \xi_2^5 = \xi_3^5 = \xi_1\xi_2\xi_3 = 1\}$
acts diagonally on $\cz^3$. One choice of resolution is the following:

Blow up the 10 intersection points of the singular lines. This produces three
exceptional divisors for each point and 30 lines of $A_1$ singularities where
two of these divisors intersect.

Blow up the 10 lines of $A_4$ singularities and the 30 lines of $A_1$ singularities.
This produces $50 = 2\cdot 10 + 1\cdot 30$ new exceptional divisors.

Blow up the remaining singular curves (intersection of divisors coming from
the blowup of the lines of $A_4$ singularities). This produces $20 = 2\cdot 10$
new exceptional divisors.

Now the singular locus consists of $60=6\cdot 10$ nodes (2 nodes on each exceptional
divisor from the first step) which can be resolved by a projective small resolution
(they are contained in the strict transforms of some of the 5 planes given by
$x_i = x_j + x_k + x_l + x_m = 0$, $\{i,j,k,l,m\}=\{0,1,2,3,4\}$).

Denote by $\tilde{Y}_{\mu}$ such a resolution of all the singularities of $Y_{\mu}$.
The Euler characteristic and Hodge numbers of $\tilde{Y}_{\mu}$ are
\[
\chi(\tilde{Y}_{\mu}) = \chi(Y_{\mu}) + 200 = 200, \quad
h^{2,1}(\tilde{Y}_{\mu}) = 1, \quad h^{1,1}(\tilde{Y}_{\mu}) = 100.
\]

On $Y=Y_1$ there is the additional node $(1:1:1:1:1)$. It is the image of the
125 nodes of the Schoen quintic $X$ under the map $\phi$.
On $X$ the node $(1:1:1:1:1)$ is contained in the smooth quadric surface
$Q$; thus on $Y$ the node $(1:1:1:1:1)$ is contained in the smooth surface
$\phi(Q)$ (which does not meet the singular lines) so there exist projective
small resolutions. Let $\tilde{Y}$ denote such a small
resolution of (all the singularities of) $Y$. The Euler characteristic of $\tilde{Y}$ is
\[
\chi(\tilde{Y}) = \chi(Y) + 200 + 1 = 202.
\]

The Hodge numbers of $\tilde{Y}$ are
\[
h^{1,1}(\tilde{Y}) = 101, \qquad h^{2,1}(\tilde{Y}) = 0.
\]

Investigating the arithmetic of $\tilde{Y}$ in detail (cf. \cite{Meyer}) we find
\begin{align*}
& \tr(\Frob_p^*|H^{3}_{\text{\'et}}(\tilde{Y},\qz_{\ell})) = a_p\\
& \qquad = \begin{cases}
p^3 + p^2 \hspace*{24,4pt} + 1 - \#Y_p, & p \equiv 1,4 \mod 5,\\
p^3 + p^2 + 2p + 1 - \#Y_p, & p \equiv 2,3 \mod 5.\\
\end{cases}
\end{align*}

The same weight four newform as for $\tilde{X}$ has to occur in the cohomology
of $\tilde{Y}$ because the map $\phi$ defines a correspondence between the two
varieties.

\section{Generalizations}

The quintic is one of the 13 hypergeometric threefolds with $h^{1,1}=1$ listed in
section 6.4 of \cite{BatvS}. It is tempting to ask if similarly nice equations
can be produced for the mirrors of some of the other cases. We will investigate
one example, namely the complete intersection of two cubics in $\pz^5$.

Let the threefold $V_{\lambda} \subset \pz^5$ be given by the equations
\begin{align*}
x_0^3 + x_1^3 + x_2^3 & = 3 \lambda x_3 x_4 x_5,\\
x_3^3 + x_4^3 + x_5^3 & = 3 \lambda x_0 x_1 x_2.
\end{align*}

This is a complete intersection which is invariant under the group
$\tilde{G}$ (of order 81) of transformations
$g_{\alpha, \beta, \delta, \epsilon, \mu}$ where $\alpha, \beta,
\delta, \epsilon \in \zz/3\zz$, $\mu \in \zz/9\zz$, and
$\mu \equiv \alpha + \beta \equiv \delta + \epsilon \mod 3$ (Note
the misprint in \cite{LibTei}). These transformations act as
\begin{align*}
g_{\alpha, \beta, \delta, \epsilon, \mu} : \; & (x_0:x_1:x_2:x_3:x_4:x_5)\\
\mapsto & (
\xi_3^{\alpha} \xi_9^{\mu} x_0:
\xi_3^{\beta} \xi_9^{\mu} x_1:
\xi_9^{\mu} x_2:
\xi_3^{-\delta} \xi_9^{-\mu} x_3:
\xi_3^{-\epsilon} \xi_9^{-\mu} x_4:
\xi_9^{-\mu} x_5)
\end{align*}
where $\xi_i$ is a fixed primitive $i$-th root of unity.
For generic $\lambda$ the variety $V_{\lambda}$ is a smooth Calabi--Yau threefold with
Euler characteristic $\chi(V_{\lambda}) = -144$. Libgober and Teitelbaum (\cite{LibTei})
show that the mirror partner of $V_{\lambda}$ can be described as a resolution of the
quotient $V_{\lambda}/\tilde{G}$.

In analogy to the case of the quintic we can study the complete intersection
$W_{\lambda}\subset\pz^5$ given by the equations
\begin{align*}
(x_0 + x_1 + x_2)^3 & = (3\lambda)^3 x_3 x_4 x_5,\\
(x_3 + x_4 + x_5)^3 & = (3\lambda)^3 x_0 x_1 x_2,
\end{align*}
but in this case the restriction of the map
\[
\psi : \pz^5 \lra \pz^5, \quad (x_0:x_1:x_2:x_3:x_4:x_5) \mapsto (x_0^3:x_1^3:x_2^3:x_3^3:x_4^3:x_5^3)
\]
to $V_{\lambda}\lra W_{\lambda}$ does not divide out the whole group $\tilde{G}$ but only
a subgroup $\tilde{H}$ with 27 elements so that (a resolution of) $W_{\lambda}$ is not the
mirror of $V_{\lambda}$. The action of $\tilde{G}/\tilde{H}$ on $W_{\lambda}$ is
multiplication of $x_0$, $x_1$, $x_2$ by $\xi_3$.

Consider the coordinate change
\begin{align*}
& x_0\mapsto x_0+x_1+x_2,\quad x_1\mapsto x_0+\xi_3 x_1+\xi_3^2 x_2,\quad x_2\mapsto x_0+\xi_3^2 x_1+\xi_3 x_2,\\
& x_3\mapsto x_3+x_4+x_5,\quad x_4\mapsto x_3+\xi_3 x_4+\xi_3^2 x_5,\quad x_5\mapsto x_3+\xi_3^2 x_4+\xi_3 x_5.
\end{align*}
In these new coordinates $W_{\lambda}$ is given by the equations
\begin{align*}
x_0^3 & = \lambda^3(x_3^3 + x_4^3 + x_5^3 - 3x_3x_4x_5),\\
x_3^3 & = \lambda^3(x_0^3 + x_1^3 + x_2^3 - 3x_0x_1x_2),
\end{align*}
or, after setting $\nu = 1/\lambda^3$ for convenience, by the equations
\begin{align*}
x_3^3 + x_4^3 + x_5^3 - \nu x_0^3 & = 3 x_3 x_4 x_5,\\
x_0^3 + x_1^3 + x_2^3 - \nu x_3^3 & = 3 x_0 x_1 x_2.
\end{align*}

Applying the map $\psi$ again we obtain the quotient $\tilde{W}_{\nu}$ of $W_{\lambda}$ given by the equations
\begin{align*}
(x_3 + x_4 + x_5 - \nu x_0)^3 & = 3^3 x_3 x_4 x_5,\\
(x_0 + x_1 + x_2 - \nu x_3)^3 & = 3^3 x_0 x_1 x_2.
\end{align*}
Unfortunately this way we divide out more than the group $\tilde{G}/\tilde{H}$ so that
(a resolution of) $\tilde{W}_{\nu}$ is a quotient of the mirror of $V_{\lambda}$.

\begin{Question}
Given a mirror pair of Calabi--Yau threefolds such that one of them can be described
as a complete intersection in a certain weighted projective space. Under which conditions
is there a description of the mirror as (a resolution of) a complete intersection of the same type?
\end{Question}

\end{document}